\documentclass[12pt,leqno]{article}
\usepackage{latexsym,amsfonts,amsmath,amssymb}
\newtheorem{theorem}{Theorem}
\newtheorem{corollary}[theorem]{Corollary}
\newtheorem{sublemma}{Lemma}[theorem]
\newtheorem{lemma}[theorem]{Lemma}
\newtheorem{question}[theorem]{Question}
\newtheorem{observation}[theorem]{Observation}
\newtheorem{claim}[theorem]{Claim}
\newtheorem{subclaim}{Claim}[sublemma]
\newtheorem{conjecture}[theorem]{Conjecture}
\newtheorem{fact}[theorem]{Fact}
\newtheorem{definition}[theorem]{Definition}
\newtheorem{remark}[theorem]{Remark}
\newtheorem{example}[theorem]{Example}
\newtheorem{exercise}{Exercise}[section]
\def\Theorem #1.#2 #3\par{\setbox1=\hbox{#1}\ifdim\wd1=0pt
   \begin{theorem}{\rm #2} #3\end{theorem}\else
   \newtheorem{#1}[theorem]{#1}\begin{#1}\label{#1}{\rm #2} #3\end{#1}\fi}
\def\Corollary #1.#2 #3\par{\setbox1=\hbox{#1}\ifdim\wd1=0pt
   \begin{corollary}{\rm #2} #3\end{corollary}\else
   \newtheorem{#1}[theorem]{#1}\begin{#1}\label{#1}{\rm #2} #3\end{#1}\fi}
\def\Lemma #1.#2 #3\par{\setbox1=\hbox{#1}\ifdim\wd1=0pt
   \begin{lemma}{\rm #2} #3\end{lemma}\else
   \newtheorem{#1}[theorem]{#1}\begin{#1}\label{#1}{\rm #2} #3\end{#1}\fi}
\def\SubLemma #1.#2 #3\par{\setbox1=\hbox{#1}\ifdim\wd1=0pt
   \begin{sublemma}{\rm #2} #3\end{sublemma}\else
   \newtheorem{#1}{#1}[theorem]\begin{#1}\label{#1}{\rm #2} #3\end{#1}\fi}
\def\Question #1.#2 #3\par{\setbox1=\hbox{#1}\ifdim\wd1=0pt
   \begin{question}{\rm #2} #3\end{question}\else
   \newtheorem{#1}[theorem]{#1}\begin{#1}\label{#1}{\rm #2} #3\end{#1}\fi}
\def\Observation #1.#2 #3\par{\setbox1=\hbox{#1}\ifdim\wd1=0pt
   \begin{observation}{\rm #2} #3\end{observation}\else
   \newtheorem{#1}[theorem]{#1}\begin{#1}\label{#1}{\rm #2} #3\end{#1}\fi}
\def\Claim #1.#2 #3\par{\setbox1=\hbox{#1}\ifdim\wd1=0pt
   \begin{claim}{\rm #2} #3\end{claim}\else
   \newtheorem{#1}[theorem]{#1}\begin{#1}\label{#1}{\rm #2} #3\end{#1}\fi}
\def\SubClaim #1.#2 #3\par{\setbox1=\hbox{#1}\ifdim\wd1=0pt
   \begin{subclaim}{\rm #2} #3\end{subclaim}\else
   \newtheorem{#1}{#1}[sublemma]\begin{#1}\label{#1}{\rm #2} #3\end{#1}\fi}
\def\Conjecture #1.#2 #3\par{\setbox1=\hbox{#1}\ifdim\wd1=0pt
   \begin{conjecture}{\rm #2} #3\end{conjecture}\else
   \newtheorem{#1}[theorem]{#1}\begin{#1}\label{#1}{\rm #2} #3\end{#1}\fi}
\def\Fact #1.#2 #3\par{\setbox1=\hbox{#1}\ifdim\wd1=0pt
   \begin{fact}{\rm #2} #3\end{fact}\else
   \newtheorem{#1}[theorem]{#1}\begin{#1}\label{#1}{\rm #2} #3\end{#1}\fi}
\def\Definition #1.#2 #3\par{\setbox1=\hbox{#1}\ifdim\wd1=0pt
   \begin{definition}{\rm #2} {\rm #3}\end{definition}\else
   \newtheorem{#1}[theorem]{#1}\begin{#1}\label{#1}{\rm #2} {\rm #3}\end{#1}\fi}
\def\Remark #1.#2 #3\par{\setbox1=\hbox{#1}\ifdim\wd1=0pt
   \begin{remark}{\rm #2} {\rm #3}\end{remark}\else
   \newtheorem{#1}[theorem]{#1}\begin{#1}\label{#1}{\rm #2} {\rm #3}\end{#1}\fi}
\def\Example #1.#2 #3\par{\setbox1=\hbox{#1}\ifdim\wd1=0pt
   \begin{example}{\rm #2} #3\end{example}\else
   \newtheorem{#1}[theorem]{#1}\begin{#1}\label{#1}{\rm #2} #3\end{#1}\fi}
\def\Exercise #1.#2 #3\par{\setbox1=\hbox{#1}\ifdim\wd1=0pt
   {\footnotesize\begin{exercise}{\rm #2} {\rm #3}\end{exercise}}\else
   \newtheorem{#1}[section]{#1}{\footnotesize\begin{#1}\label{#1}{\rm #2} {\rm #3}\end{#1}}\fi}
\def\QuietTheorem #1.#2 #3\par{\setbox1=\hbox{#1}\ifdim\wd1=0pt\proclaim{Theorem {\rm #2}}{#3}\else\proclaim{#1 {\rm #2}}{#3}\fi}
\newcommand{\proclaim}[2]{\smallskip\noindent{\bf #1} {\sl#2}\par\smallskip}
\def\Proclaim #1.#2 #3\par{\proclaim{#1 {\rm #2}}{#3}}
\newenvironment{proof}{\noindent}{\kern2pt\QEDbox\par\bigskip}
\def\Proof#1: {\setbox1=\hbox{#1}\ifdim\wd1=0pt\begin{proof}{\bf Proof: }\else\medskip\begin{proof}{\bf #1: }\fi}
\newcommand{\QED}{\end{proof}}
\def\BF#1.{{\bf #1.}}
\def\Abstract #1\par{\begin{quotation}{\singlespaced\footnotesize{\noindent{\bf Abstract.~}#1}}\end{quotation}}
\def\Title #1\par{\title{#1}\maketitle}
\def\Author #1\par{\author{#1}}
\def\Acknowledgement#1\par{\thanks{#1}}
\def\Chapter #1\par{\chapter{#1}}
\def\Section #1\par{\section{#1}}
\def\QuietSection #1\par{\section*{#1}}
\def\SubSection #1\par{\subsection{#1}}
\def\SubSubSection #1\par{\subsubsection{#1}}
\def\MidTitle #1\par{\bigskip\goodbreak\centerline{\small\bf #1}\bigskip\noindent}
\def\Margin #1\par{\marginpar{\tiny #1}}

\newcommand{\singlespaced}{\baselineskip=15pt}
\def\bottomnote #1\par{{\renewcommand{\thefootnote}{}\footnotetext{#1}}}
\renewcommand{\P}{{\mathbb P}}
\newcommand{\Q}{{\mathbb Q}}

\newcommand{\Mbar}{{\overline{M}}}
\newcommand{\Nbar}{{\overline{N}}}
\newcommand{\Vbar}{{\overline{V}}}
\newcommand{\Xbar}{{\overline{X}}}

\newcommand{\Qdot}{{\dot\Q}}
\newcommand{\qdot}{{\dot q}}

\newcommand{\rdot}{{\dot r}}
\newcommand{\sdot}{{\dot s}}

\newcommand{\id}{\mathop{\hbox{\small id}}}

\newfont{\msam}{msam10 at 12pt}

\newcommand{\of}{\subseteq}
\newcommand{\ofnoteq}{\subsetneq}

\newcommand{\set}[1]{\{\,{#1}\,\}}

\newcommand{\elesub}{\prec}

\newcommand{\cof}{\mathop{\rm cof}}

\newcommand{\image}{\mathbin{\hbox{\tt\char'42}}}
\newcommand{\plus}{{+}}

\newcommand{\restrict}{\upharpoonright}
\newcommand{\satisfies}{\models}
\newcommand{\forces}{\Vdash}

\newcommand{\concat}{\mathbin{{}^\smallfrown}}

\newcommand{\union}{\cup}
\newcommand{\Union}{\bigcup}
\newcommand{\intersect}{\cap}

\newcommand{\Pforces}{\forces_{\P}}

\newcommand{\trianglelt}{\lhd}

\newcommand{\smalllt}{\mathrel{\mathchoice{\raise2pt\hbox{$\scriptstyle<$}}{\raise1pt\hbox{$\scriptstyle<$}}{\scriptscriptstyle<}{\scriptscriptstyle<}}}
\newcommand{\smallleq}{\mathrel{\mathchoice{\raise2pt\hbox{$\scriptstyle\leq$}}{\raise1pt\hbox{$\scriptstyle\leq$}}{\scriptscriptstyle\leq}{\scriptscriptstyle\leq}}}
\newcommand{\ltomega}{{{\smalllt}\omega}}

\newcommand{\ltkappa}{{{\smalllt}\kappa}}

\newcommand{\lttheta}{{{\smalllt}\theta}}

\newcommand{\leqdelta}{{{\smallleq}\delta}}
\newcommand{\ltdelta}{{{\smalllt}\delta}}
\newcommand{\ltbeta}{{{\smalllt}\beta}}

\newcommand{\card}[1]{{|#1|}}
\newcommand{\boolval}[1]{\mathopen{\lbrack\!\lbrack}\,#1\,\mathclose{\rbrack\!\rbrack}}
\def\[#1]{\boolval{#1}}

\newcommand{\UnderTilde}[1]{{\setbox1=\hbox{$#1$}\baselineskip=0pt\vtop{\hbox{$#1$}\hbox to\wd1{\hfil$\sim$\hfil}}}{}}
\newcommand{\Undertilde}[1]{{\setbox1=\hbox{$#1$}\baselineskip=0pt\vtop{\hbox{$#1$}\hbox to\wd1{\hfil$\scriptstyle\sim$\hfil}}}{}}
\newcommand{\undertilde}[1]{{\setbox1=\hbox{$#1$}\baselineskip=0pt\vtop{\hbox{$#1$}\hbox to\wd1{\hfil$\scriptscriptstyle\sim$\hfil}}}{}}
\newcommand{\UnderdTilde}[1]{{\setbox1=\hbox{$#1$}\baselineskip=0pt\vtop{\hbox{$#1$}\hbox to\wd1{\hfil$\approx$\hfil}}}{}}
\newcommand{\Underdtilde}[1]{{\setbox1=\hbox{$#1$}\baselineskip=0pt\vtop{\hbox{$#1$}\hbox to\wd1{\hfil\scriptsize$\approx$\hfil}}}{}}

\newcommand{\st}{\mid}
\renewcommand{\th}{{\hbox{\scriptsize th}}}

\newcommand{\minus}{\setminus}
\newcommand{\iso}{\cong}
\def\<#1>{\langle#1\rangle}
\newcommand{\ot}{\mathop{\rm ot}\nolimits}
\newcommand{\QEDbox}{\fbox{}}
\newcommand{\cp}{\mathop{\rm cp}}

\newcommand{\ORD}{\mathop{\hbox{\sc ord}}}

\newcommand{\ZFC}{\hbox{\sc zfc}}

\newcommand{\cell}[1]{\boxit{\hbox to 17pt{\strut\hfil$#1$\hfil}}}
\newcommand{\head}[2]{\lower2pt\vbox{\hbox{\strut\footnotesize\it\hskip3pt#2}\boxit{\cell#1}}}
\newcommand{\boxit}[1]{\setbox4=\hbox{\kern2pt#1\kern2pt}\hbox{\vrule\vbox{\hrule\kern2pt\box4\kern2pt\hrule}\vrule}}
\newcommand{\Col}[3]{\hbox{\vbox{\baselineskip=0pt\parskip=0pt\cell#1\cell#2\cell#3}}}
\newcommand{\tapenames}{\raise 5pt\vbox to .7in{\hbox to .8in{\it\hfill input: \strut}\vfill\hbox to
.8in{\it\hfill scratch: \strut}\vfill\hbox to .8in{\it\hfill output: \strut}}}
\newcommand{\Head}[4]{\lower2pt\vbox{\hbox to25pt{\strut\footnotesize\it\hfill#4\hfill}\boxit{\Col#1#2#3}}}
\newcommand{\Dots}{\raise 5pt\vbox to .7in{\hbox{\ $\cdots$\strut}\vfill\hbox{\ $\cdots$\strut}\vfill\hbox{\
$\cdots$\strut}}}
\renewcommand{\dots}{\raise5pt\hbox{\ $\cdots$}}
\newcommand{\factordiagramup}[6]{$$\begin{array}{ccc}
#1&\raise3pt\vbox{\hbox to60pt{\hfill$\scriptstyle
#2$\hfill}\vskip-6pt\hbox{$\vector(4,0){60}$}}&#3\\ \vbox
to30pt{}&\raise22pt\vtop{\hbox{$\vector(4,-3){60}$}\vskip-22pt\hbox
to60pt{\hfill$\scriptstyle #4\qquad$\hfill}}
     &\ \ \lower22pt\hbox{$\vector(0,3){45}$}\ {\scriptstyle #5}\\
\vbox to15pt{}&&#6\\
\end{array}$$}
\newcommand{\factordiagram}[6]{$$\begin{array}{ccc}
#1&&\\ \ \ \raise22pt\hbox{$\vector(0,-3){45}$}\ {\scriptstyle #2}
&\raise22pt\hbox{$\vector(2,-1){90}$}\raise5pt\llap{$\scriptstyle#3$\qquad\quad}&\vbox
to25pt{}\\ #4&\raise3pt\vbox{\hbox to90pt{\hfill$\scriptstyle
#5$\hfill}\vskip-6pt\hbox{$\vector(4,0){90}$}}&#6\\
\end{array}$$}
\newcommand{\df}{\it} 
\begin{document}
\author{Joel David Hamkins\\
\normalsize\sc The City University of New York}
\date{}
\bottomnote MSC: 03E55, 03E40. Keywords: Large cardinals, forcing. I am affiliated with the College of Staten Island of CUNY and with The CUNY
Graduate Center. My research has been supported by grants from Georgia State University, the Research Foundation of CUNY and the National Science
Foundation. I would like to thank Arthur Apter for his patient encouragement and the referee for some insightful remarks.

\Title Extensions with the approximation and cover properties have no new large cardinals

\Abstract If an extension $V\of\Vbar$ satisfies the $\delta$ approximation and cover properties for classes and $V$ is a class in $\Vbar$, then every
suitably closed embedding $j:\Vbar\to\Nbar$ in $\Vbar$ with critical point above $\delta$ restricts to an embedding $j\restrict V$ amenable to the
ground model $V$. In such extensions, therefore, there are no new large cardinals above $\delta$. This result extends work in
\cite{Hamkins2001:GapForcing}.

\Section Introduction

While an important theme in set theory concerns the preservation of large cardinals from a ground model to various forcing extensions, set theorists
often expect conversely that a forcing extension will not exhibit new instances of large cardinals. After all, the smallest large cardinals, such as
inaccessible and Mahlo cardinals, are downwards absolute to any model; those large cardinals not implying $0^\sharp$ are downwards absolute to $L$,
and many stronger notions are downwards absolute to the core models. Kunen \cite{Kunen78:SaturatedIdeals} discovered, however, that forcing sometimes
can create new large cardinals: a non-weakly compact cardinal $\kappa$ can become measurable or more after adding a branch to a $\kappa$ Suslin tree.
Other examples show that adding even a Cohen subset to a non-measurable cardinal $\kappa$ can make it supercompact or more.

Despite these examples, the Main Theorem of this article confirms the general expectation by showing that for a large class of extensions
$V\of\Vbar$, every suitably closed embedding $j:\Vbar\to\Nbar$ in the extension $\Vbar$ lifts an embedding $j\restrict V:V\to N$ amenable to the
ground model. Since these ground model embeddings typically witness the corresponding large cardinal property in $V$, it follows that the extension
$\Vbar$ has no new large cardinals. This work generalizes \cite{Hamkins2001:GapForcing}.

The Main Theorem therefore concerns the {\df lifting property} for $V\of\Vbar$, the property asserting that every suitable embedding in $\Vbar$ lifts
an embedding in $V$. This property is of course already well known in many cases, often with little or no restriction on $\Vbar$. For example,
$0^\sharp$ cannot be added by set forcing over $L$, and every embedding $j:L[\mu]\to L[j(\mu)]$ in a forcing extension of $L[\mu]$ is necessarily an
iteration of $\mu$. Similar results hold for larger cardinals with respect to the core models.

Nevertheless, there are easy counterexamples to the lifting property when the embeddings lack closure. If there are two normal measures on a
measurable cardinal $\kappa$, for example, and $x$ is a Cohen real, then in $V[x]$ the iteration $j$ of the extensions of these measures, chosen by
the digits of $x$, cannot lift an embedding amenable to $V$, because from $j\restrict P(\kappa)^V$ one easily reconstructs $x$. Jensen observed (in
the 1980s) that the lifting property can fail without the closure requirement, by pointing out that the core model is not the union of all mice when
there are mice with more than one normal measure.

Other counterexamples satisfy the closure requirement. For example, many large cardinals $\kappa$ are preserved by the forcing to add a Cohen subset
$A\of\kappa$; but no embedding $j:V[A]\to N[j(A)]$ can lift an embedding amenable to $V$, because necessarily $A=j(A)\intersect\kappa\in N[j(A)]$ and
so $A\in N$ by the closure of $j(A)$, leading to $N\not\of V$. By adding a Cohen subset to every inaccessible cardinal below $\kappa$, and then
finally at $\kappa$, one can arrange that a large cardinal $\kappa$ is killed in $V[G]$ and resurrected in $V[G][A]$, leading to a strong violation
of the lifting property for $V[G]\of V[G][A]$.

Several open questions remain. The extent of the lifting property for extensions without the approximation and cover properties is not known.
Countably (strategically) closed forcing, such as many forward Easton products, can usually be handled by first adding a Cohen real, so that the
combined forcing has a closure point at $\omega$.  For forcing without a closure point, the question is largely open. Gitik has constructed a
counterexample to the lifting property for iterated Prikry forcing, perhaps a worst-case example for lacking closure points. It is open whether one
can reduce the closure requirement of the Main Theorem to $\Nbar^\ltdelta\of\Nbar$, which would allow for $\Nbar^\omega\of\Nbar$ in extensions with
the countable approximation and cover properties.

\Section The Main Theorem

Let me now state and prove the main theorem. By a model of set theory, I mean a model of some fixed large finite fragment of \ZFC, sufficiently
powerful to carry out such standard arguments as the construction of the cumulative hierarchy $V_\alpha$, Mostowski collapses and so on. For
definiteness, take it to mean a model of the $\Sigma_{100}$ fragment of \ZFC. Throughout this article, $V\of\Vbar$ is an extension consisting of two
transitive class models of \ZFC, viewed as the respective universes of all sets, with the principal example occurring when $\Vbar$ is a forcing
extension of $V$.

\Definition. A pair of transitive classes $M\of N$ satisfies the {\df $\delta$ approximation property} if whenever $A\of M$ is a set in $N$ and
$A\intersect a\in M$ for any $a\in M$ of size less than $\delta$ in $M$, then $A\in M$. For models of set theory equipped with classes, the pair
$M\of N$ satisfies the $\delta$ approximation property {\df for classes} if whenever $A\of M$ is a class of $N$ and $A\intersect a\in M$ for any $a$
of size less than $\delta$ in $M$, then $A$ is a class of $M$. I will refer to the sets $A\intersect a$, where $a$ has size less than $\delta$ in
$M$, as the {\df $\delta$ approximations} to $A$ over $M$.

\Definition. The pair $M\of N$ satisfies the {\df $\delta$ cover property} if for every set $A$ in $N$ with $A\of M$ and $|A|^N<\delta$, there is a
set $B\in M$ with $A\of B$ and $|B|^M<\delta$.

\Theorem Main Theorem. Suppose that $V\of \Vbar$ satisfies the $\delta$ approximation and cover properties, $\delta$ is regular, $\Mbar$ is a
transitive submodel of $\Vbar$ such that $M=\Mbar\intersect V$ is also a model of set theory, and $j:\Mbar\to \Nbar$ is a (possibly external) cofinal
elementary embedding of $\Mbar$ into a transitive class $\Nbar\of\Vbar$. Suppose further that $\delta<\cp(j)$, $P(\delta)^\Vbar\of \Mbar$ and that
$\Mbar^\ltdelta\of\Mbar$ and $\Nbar^\delta\of\Nbar$ in $\Vbar$. Let $N=\Union j\image M$, so that $j\restrict M:M\to N$. Then:
\begin{enumerate}
 \item If\/ $\Mbar$ is a set in $\Vbar$, then $M$ is a set in $V$.
 \item $N\of V$; indeed, $N=\Nbar\intersect V$.
 \item If $j$ is amenable to $\Vbar$,
 then $j\restrict M$ is amenable to $V$. In particular, if $j$ is a set in $\Vbar$, then the restricted embedding $j\restrict M$ is a set in $V$.
 \item If $j$ and $M$ are classes in $\Vbar$ and $V\of\Vbar$ satisfies the $\delta$ approximation property for classes, then $j\restrict M$ is a
  class of $V$.
\end{enumerate}

\Proof: Let me focus at first on the central case, where $\Mbar=\Vbar$ and consequently $M=V$. After this, I'll explain how to modify the argument
for the general case. In the central case, we have an embedding $j:\Vbar\to\Nbar$ with $\Nbar\of\Vbar$ and $\Nbar^\delta\of\Nbar$ in $\Vbar$. If
$N=\Union j\image V$, then a simple induction on formulas shows that the restricted embedding $j\restrict V:V\to N$ is elementary.

\SubLemma. $N\of \Nbar$ satisfies the $\delta$ approximation and cover properties.\label{NNbar}

\Proof: We apply the elementarity of $j$ to the corresponding properties for $V\of\Vbar$. Specifically, suppose that $A\in\Nbar$ and all $\delta$
approximations to $A$ over $N$ are in $N$. Since the embedding is cofinal, there is some $B\in\Vbar$ with $A\of j(B)$. In $\Vbar$, any subset of $B$
having all $\delta$ approximations over $P(B)^V$ in $P(B)^V$, is in $P(B)^V$. Thus, the corresponding fact is true in $\Nbar$ about subsets of $j(B)$
and their $j(\delta)=\delta$ approximations over $j(P(B)^V)=P(j(B))^N$. In particular, since $A\of j(B)$ has all its $\delta$ approximation over $N$
in $N$, I conclude that $A\in N$, as desired. The $\delta$ cover property is similar.\QED


\SubLemma. If $A\of \ORD^N$ is a set of size less than $\delta$ in $\Vbar$, then there is a set $B\in V\intersect N$ of size at most $\delta$ with
$A\of B$.\label{VNcover}

\Proof: Suppose that $A=A_0\of\ORD^N$ has size less than $\delta$ in $\Vbar$. It follows that $A\in\Nbar$ and so by Lemma \ref{NNbar} there is a set
of ordinals $A_1\in N$ of size less than $\delta$ with $A_0\of A_1$. Since also $A_1\in \Vbar$ there is a set $A_2\in V$ of size less than $\delta$
with $A_1\of A_2$. We may continue bouncing between $N\of \Nbar$ and $V\of\Vbar$, using the regularity of $\delta$ at limit stages, in order to build
a sequence $\<A_\alpha\st \alpha<\delta>$ in $\Vbar$ such that $\alpha<\beta\implies A_\alpha\of A_\beta$, all $A_\alpha$ are subsets of $\ORD^N$
having size less than $\delta$, and unboundedly often $A_\alpha\in V$ and unboundedly often $A_\alpha\in N$. Let $B=\Union_{\alpha<\delta}A_\alpha$.
Since $\Nbar^\delta\of\Nbar$ and $B$ has size at most $\delta$ in $\Vbar$, we conclude that $B$ is in $\Nbar$ and has size at most $\delta$ there. If
$a$ is any set of ordinals of size less than $\delta$ in $V$, then $B\intersect a=A_\alpha\intersect a$ for sufficiently large $\alpha$, and so
$B\intersect a\in V$. Thus, all the $\delta$ approximations to $B$ over $V$ are in $V$, and so $B\in V$. Similarly, all the $\delta$ approximations
to $B$ over $N$ are in $N$, and so $B\in N$. Therefore $B\in V\intersect N$, as desired.\QED

\SubLemma. $V$ and $N$ have the same subsets of $\ORD^N$ of size less than $\delta$.\label{VNdelta}

\Proof: Suppose that $A\of\ORD^N$ has size less than $\delta$ in $\Vbar$. By Lemma \ref{VNcover} there is a set $B\in V\intersect N$ of size at most
$\delta$ in $\Vbar$ with $A\of B$. Enumerate $B=\set{\beta_\alpha\st\alpha<\bar\delta}$ in the natural order, where $\bar\delta=\ot(B)<\delta^\plus$,
and let $a=\set{\alpha<\bar\delta\st \beta_\alpha\in A}$. If $A$ is in either $V$ or $N$, then so is $a$, since it is constructible from $A$ and $B$.
Since $\bar\delta$ is below the critical point of $j$, we know that $j(a)=a$. Since $j(P(\bar\delta)^V)=P(\bar\delta)^N$, it follows that $a\in V$ if
and only if $a=j(a)\in N$. So $a$ must be in both $V$ and $N$. Finally, as $A$ is constructible from $B$ and $a$, we conclude that $A$ is in both $V$
and $N$ as well.\QED

\SubLemma. $N\of V$. Indeed, $N=\Nbar\intersect V$.\label{N=NbarIntersectV}

\Proof: For the forward inclusion, it suffices to show that every set of ordinals in $N$ is in $V$. Suppose that $A\of \ORD^N$ and $A\in N$. Fix any
$a\in V$ of size less than $\delta$ in $V$, and consider $A\intersect a$. We may assume $a\of\ORD^N$. It follows by Lemma \ref{VNdelta} that $a\in N$
and so also $A\intersect a\in N$. By Lemma \ref{VNdelta} again, it follows that $A\intersect a\in V$, and so every $\delta$ approximation to $A$ over
$V$ is in $V$. Consequently, by the $\delta$ approximation property, $A\in V$, as desired.

Conversely, suppose that $A\in\Nbar\intersect V$, considering first the case when $A$ is a set of ordinals. If $a\of\ORD$ has size less than $\delta$
in $N$, then $a\in V$ by Lemma \ref{VNdelta}, and so $A\intersect a\in V$. Thus, $A\intersect a\in N$ by Lemma \ref{VNdelta} again, and so all the
$\delta$ approximations to $A$ over $N$ are in $N$. By the $\delta$ approximation property of $N\of\Nbar$, we conclude $A\in N$, as desired. For the
general case, suppose that $A$ is any set in $\Nbar\intersect V$. By $\in$-induction, suppose that every element of $A$ is in $N$. Thus, $A\of B$ for
some set $B\in N$. Enumerate $B=\set{b_\alpha\st \alpha<\beta}$ in $N\of V$, and consider $A_0=\set{\alpha<\beta\st b_\alpha\in A}$. This is
constructible from $A$ and the enumeration of $B$, and so it is in both $\Nbar$ and $V$. Therefore, $A_0\in N$ by the earlier argument of this
paragraph. And since $A$ is constructible from $A_0$ and the enumeration of $B$, we conclude $A\in N$, as desired.\QED

\SubLemma. If $j$ is amenable to $\Vbar$, then $j\restrict V$ is amenable to $V$.\label{amenable}

\Proof: Assume that $j$ is amenable to $\Vbar$ and suppose $A\in V$. In order to show $j\restrict A\in V$, it suffices to show that all $\delta$
approximations to $j\restrict A$ over $V$ are in $V$. And for this, it suffices to show that $j\restrict a\in V$ for any $a$ of size less than
$\delta$ in $V$. Enumerate $a$ as $\vec a=\<a_\alpha\st\alpha<\beta>$ in $V$, and observe that $j(\vec a)=\<j(a_\alpha)\st\alpha<\beta>$ because
$\beta<\delta<\cp(j)$. Since $j(\vec a)\in N\of V$, we may construct $j\restrict a=\set{\<a_\alpha,j(a_\alpha)>\st\alpha<\beta}$ from $\vec a$ and
$j(\vec a)$ in $V$.\QED

In particular, if $j$ and $V$ are classes in $\Vbar$ and we have the $\delta$ approximation property for classes, then the previous argument shows
that all $\delta$ approximations to $j\restrict V$ over $V$ are in $V$, and so $j\restrict V$ is a class in $V$. This completes the proof of the Main
Theorem for the special case when $\Mbar=\Vbar$ and $M=V$.

Let me now describe the modifications that are required for the general proof. Let $N=\Union j\image M$, so that $j\restrict M:M\to N$ is an
elementary embedding.

\SubLemma. $M\of\Mbar$ satisfies the $\delta$ approximation and cover properties.\label{MMbar}

\Proof: For the $\delta$ approximation property, suppose that $A\in\Mbar$, $A\of M$ and $A\intersect a\in M$ whenever $a$ has size less than $\delta$
in $M$. Fix any $\sigma$ of size less than $\delta$ in $V$, and let $a=\sigma\intersect M$. Since this is the same as $\sigma\intersect(V_\beta)^M$
for sufficiently large $\beta$, we know $a\in V$. Since $a\of M$ has size less than $\delta$, it is in $\Mbar$ and hence in $\Mbar\intersect V=M$. So
$A\intersect a\in M$. Since $A\intersect\sigma=A\intersect a$ and $M\of V$, this means that all $\delta$ approximations to $A$ over $V$ are in $V$,
and so $A\in V$ by the $\delta$ approximation property of $V\of\Vbar$. Thus, $A\in\Mbar\intersect V=M$, as desired.

For the $\delta$ cover property, suppose that $A\of M$ has size less than $\delta$ in $\Mbar$. Since $A\of V$, $A\in\Vbar$ and $A$ has size less than
$\delta$ in $\Vbar$, there is a set $B_0$ of size less than $\delta$ in $V$ with $A\of B_0$. Using a sufficiently large $(V_\beta)^M$, there is a set
$B_1\in M\of V$ with $A\of B_1$. Thus, $A\of B_0\intersect B_1$ and $B_0\intersect B_1\of M$ has size less than $\delta$ in $V$. It follows that
$B_0\intersect B_1\in\Mbar$ and consequently in $\Mbar\intersect V=M$. Furthermore, any bijection witnessing that this set has size less than
$\delta$ in $V$ will be in $\Mbar$ and hence in $\Mbar\intersect V=M$ as well.\QED

Next, I claim that if $\Mbar$ is a set in $\Vbar$, then $M$ is a set in $V$. This is because all the $\delta$ approximations to $M$ over $V$ are in
$V$: if $a$ has size less than $\delta$ in $V$, then $M\intersect a\of B$ for some $B\in M\of V$, and so $M\intersect a=B\intersect a\in V$. In the
general case, one proves Lemma \ref{NNbar} by applying $j$ to Lemma \ref{MMbar} rather than to the inclusion $V\of\Vbar$, and in Lemma \ref{VNdelta}
one uses the hypothesis that $P(\delta)^\Vbar\of\Mbar$ in order to know that $a\in\Mbar$ and also $\bar\delta<\cp(j)$, which gives $j(a)=a$, so that
$a\in V$ if and only if $a\in N$. In Lemma \ref{amenable}, one shows that $j\restrict M$ has all its $\delta$ approximations in $V$. It follows that
if $j$ is amenable to $\Vbar$, then $j\restrict M$ is amenable to $V$. For the same reason, if $j$ and $M$ are classes in $\Vbar$ and one has the
$\delta$ approximation property for classes, then $j\restrict M$ is a class in $V$. So if $j$ is a set in $\Vbar$, then $j\restrict M$ is a set in
$V$, without any consideration of classes. This completes the proof of the Main Theorem.\QED

Let me prove a bit more about the situation of the main theorem.

\Corollary. Under the hypothesis of the theorem, for any $\lambda$,
\begin{enumerate}
 \item If $\Nbar^\lambda\of\Nbar$ in $\Vbar$, then $N^\lambda\of N$ in $V$.
 \item If $V_\lambda\of\Nbar$, then $V_\lambda\of N$.
\end{enumerate}\label{Extra}

\Proof: For 1, any $\lambda$ sequence over $N$ in $V$ is in $\Nbar\intersect V$, and hence in $N$. For 2, if $V_\lambda$ is a subset of $\Nbar$, then
it is a subset of $\Nbar\intersect V=N$.\QED

\Remark. The assumption in the Main Theorem that $\Nbar^\delta\of\Nbar$ in $\Vbar$ can be weakened to the assumption that $\Nbar\of\Vbar$ satisfies
the $\delta^\plus$ cover property, that is, that every subset of $\Nbar$ of size $\delta$ in $\Vbar$ is covered by an element of $\Nbar$ of size
$\delta$ in $\Nbar$. With the other hypotheses, this cover property is equivalent to $\Nbar^\delta\of\Nbar$, because if $\sigma\of\tau$ and $\tau$
has size $\delta$ in $\Nbar$, then one can enumerate $\tau=\set{b_\alpha\st\alpha<\delta}$ in $\Nbar$, and the set $\sigma$ is picked out by a
certain subset of $\delta$, which must be in $\Mbar$ and hence in $\Nbar$.

The central case is summarized in the following corollary.

\Corollary. Suppose that $V\of\Vbar$ satisfies the $\delta$ approximation and cover properties for classes. If $V$ is a class in $\Vbar$ and
$j:\Vbar\to\Nbar$ is a class embedding in $\Vbar$ with $\delta<\cp(j)$ and $\Nbar^\delta\of\Nbar$ in $\Vbar$, then the restriction $j\restrict V:V\to
N$, where $N=\Nbar\intersect V$, is a class elementary embedding in the ground model.\label{Simplified}

Additional simplifications are possible when $\Vbar=V[G]$ is a set forcing extension of $V$ and we equip the models with only their definable classes
(using the term {\df definable} here to mean definable from parameters).

\Lemma. Suppose that $V\of V[G]$ is a set forcing extension satisfying the $\delta$ approximation property (for sets). If the models are equipped
with only their definable classes, allowing a predicate for $V$ in $V[G]$, then $V\of V[G]$ also satisfies the $\delta$ approximation property for
classes.\label{SetClass}

\Proof: Suppose that $A\of V$ is a class in $V[G]$ all of whose $\delta$ approximations over $V$ are in $V$. For any ordinal $\eta$, let
$A_\eta=A\intersect V_\eta$. The $\delta$ approximations to $A_\eta$ over $V$ have the form $A_\eta\intersect a$ for some $a\in V$ of size less than
$\delta$ in $V$. But $A_\eta\intersect a=(A\intersect V_\eta)\intersect a=(A\intersect a)\intersect V_\eta$, which is the intersection of two sets in
$V$ and consequently in $V$. Thus, by the $\delta$ approximation property for sets, it follows that $A_\eta\in V$ for all $\eta$. Since we have
assumed that $A$ is definable in $V[G]$, there is a formula $\varphi$ (allowing a predicate for the ground model) and parameter $z$ such that
$V[G]\satisfies x\in A\iff \varphi(x,z)$. Let $\dot z$ be a name for $z$. Since $A\intersect V_\eta=A_\eta\in V$, there is some condition $p_\eta\in
G$ such that $x\in A_\eta\iff p_\eta\forces \varphi(\check x,\dot z)$. The mapping $\eta\mapsto p_\eta$ exists in $V[G]$ and $G$ is a set in $V[G]$,
so for unboundedly many $\eta$ the value of $p_\eta$ is the same. Let $p^*$ be this common value. It follows that $p^*$ could be used for any
$p_\eta$, and so we have for any $\eta$ that $x\in A_\eta\iff p^*\forces \varphi(\check x,\dot z)$. Thus, $x\in A\iff p^*\forces\varphi(\check x,\dot
z)$ provides a definition of $A$ as a class of $V$, using parameters $\dot z$, $p^*$ and the forcing poset.\QED

By Lemma \ref{SetClass}, the need to consider classes explicitly in set forcing extensions falls away, and the central case becomes the following.

\Corollary. If $V\of V[G]$ is a set forcing extension with the $\delta$ approximation and cover properties and $j:V[G]\to \Nbar$ is a definable
embedding in $V[G]$ with $\Nbar^\delta\of\Nbar$ and $\delta<\cp(j)$, then the restriction $j\restrict V:V\to N$, where $N=\Nbar\intersect V$, is an
elementary embedding definable in the ground model.

One can focus on the topic of the Main Theorem through the lens of measures and extenders rather than through the embeddings to which they give rise.
Definition \ref{LiftExtendDefinition} may help to clarify matters. The term {\df measure} here means any countably complete ultrafilter on any set;
so this includes supercompactness and hugeness measures along with ordinary measures on a measurable cardinal. For any measure $\mu$, let $j_\mu:V\to
M$ be the corresponding ultrapower embedding. If $\cp(j)=\kappa$, then it is easy to see that $M^\kappa\of M$.

\Definition. If $V\of\Vbar$ are two models of set theory with measures $\mu\in V$ and $\nu\in\Vbar$, then $\mu$ {\df lifts} to $\nu$ if
$j_\mu=j_\nu\restrict V$ and $\mu$ {\df extends} to $\nu$ if $\mu\of\nu$.\label{LiftExtendDefinition}

The {\it lift} and {\it extend} relations, though closely related, are in general independent. For two normal measures, lifting implies extending,
but the converse can fail, and forcing creates a number of interesting possibilities. If $\kappa$ is measurable, for example, then there is a forcing
extension $V\of V[f]$ where every $\kappa$-complete measure on $\kappa$ lifts to a normal measure in $V[f]$. There are other forcing extensions where
every measure in $V$ extends to $V[G]$, but none lift. If $\nu$ is a measure in $\Vbar$ concentrating on a set $D$ in $V$ and the ultrapower
$j:\Vbar\to\Nbar$ by $\nu$ lifts an embedding $j\restrict V$ amenable to $V$, then from $j\restrict P(D)^V$ in $V$ one can define a measure $\mu$ on
$D$ by $X\in \mu\iff s\in j(X)$, where $s=[\id]_\nu$, and it is easy to see that $\mu=\nu\intersect V$. So $\mu$ extends to $\nu$. Therefore, if a
measure $\nu$ in $\Vbar$ concentrates on a set in $V$ and $j_\nu\restrict V$ is amenable to $V$, then $\nu$ extends a measure in $V$.

Because of this, the Main Theorem implies the {\df measure extension property} for $V\of\Vbar$, namely, every $\delta^\plus$ complete measure $\nu$
in $\Vbar$ concentrating on any set in $V$ extends a measure $\mu$ in $V$. Theorem \ref{FilterExtension} proves this directly, generalized to include
filters and with weaker hypotheses than the Main Theorem. Counterexamples show, however, that $j_\nu\restrict V$ may not be the ultrapower by $\mu$,
or indeed the ultrapower by any measure or extender at all, even when $\nu$ is a normal measure in $\Vbar$ and $j_\nu\restrict V$ is definable in
$V$.

\Theorem. Suppose $\delta\leq\kappa$ and $V\of \Vbar$ satisfies the $\delta$ approximation property. If $\nu$ is a $\kappa$-complete filter in
$\Vbar$ on a set $D$ in $V$ and $\nu$ measures every subset of $D$ in $V$, then $\nu\intersect V$ is in $V$. That is, $\nu$ extends a measure in
$V$.\label{FilterExtension}

\Proof: This proof appeared in \cite{ApterHamkins2001:IndestructibleWC} for closure point forcing, but the approximation property is enough. It
suffices to show that every $\delta$ approximation to $\nu\intersect V$ is in $V$. So suppose $\sigma\in V$ has size less than $\delta$, and consider
$\sigma\intersect (\nu\intersect V)=\sigma\intersect\nu$. We may assume that every member of $\sigma$ is a subset of $D$. Let $\sigma^*$ be obtained
by closing $\sigma$ under complements in $D$. Since $\sigma^*\intersect\nu$ is a collection of fewer than $\delta$ many sets in the filter, it
follows by the $\kappa$-completeness of $\nu$ that $A=\intersect(\sigma^*\intersect\nu)$ is in $\nu$. Choose any $a\in A$. Observe now that if
$B\in\sigma\intersect\nu$ then $A\subseteq B$ and consequently $a\in B$. Conversely, if $a\in B$ and $B\in\sigma$ then because $a\notin D\minus B$ it
follows that $A\not\subseteq D\minus B$ and so $D\minus B\notin\nu$. By the assumption that $\nu$ measures every set in $V$, we conclude that
$B\in\nu$. Thus, we have proved for $B\in\sigma$ that $B\in \nu\iff a\in B$. So $\sigma\intersect\nu$ is precisely the set of all $B\in\sigma$ with
$a\in B$, and this is certainly in $V$. Therefore, I have proved that every $\delta$ approximation to $\nu$ over $V$ is in $V$. By the $\delta$
approximation property, it follows that $\nu\intersect V\in V$.\QED

A similar result holds for extenders, by combining techniques of the Main Theorem with ideas of \cite{HamkinsWoodin2000:SmallForcing}.

\Theorem. Suppose $\delta<\kappa$ and $V\of\Vbar$ satisfies the $\delta$ approximation and cover properties. If $E$ is an extender in $\Vbar$ whose
embedding $j:\Vbar\to\Nbar$ has $\cp(j)=\kappa$ and satisfies $\Nbar^\delta\of\Nbar$, then $E\intersect V$ is an extender in $V$.

\Proof: We suppose $E$ has the form $E=\set{\<A,s>\st s\in j(A)\And s\in[\lambda]^\ltomega}$. To show $E\intersect V\in V$, it suffices to show that
all the $\delta$ approximations to $E$ over $V$ are in $V$. Fix any set $a$ of size less than $\delta$ in $V$, and consider $E\intersect a$. Let
$\sigma$ be the set of all ordinals mentioned in the second coordinate of $a$. This is a set of ordinals in $V$ of size less than $\delta$, and
consequently it is in $\Nbar$. By Lemma \ref{VNcover}, there is a set $\tau\in V\intersect N$, where $N=\Union j\image V$, of size $\delta$ in both
$V$ and $N$ such that $\sigma\of\tau$. Let $\nu=\set{X\of V_\kappa\st \tau\in j(X)}$. Since $\tau$ is in $N$, this is a $\kappa$-complete measure on
$V_\kappa$ in $\Vbar$, and so by Theorem \ref{FilterExtension}, we know that $\mu=\nu\intersect V$ is a measure in $V$. I claim now that $E\intersect
a$ is constructible from $\mu$ and $\tau$ in $V$. Suppose that $\<A,s>\in a$, and I want to determine in $V$ whether $\<A,s>\in E$. Enumerate
$s=\<\alpha_0,\ldots,\alpha_k>$, where $\alpha_i\in\tau$. Each ordinal $\alpha_i$ is the $\beta_i^\th$ element of $\tau$ for some unique $\beta_i$.
If $f(t)$ is the finite sequence consisting of the $\beta_0^\th$, \ldots, $\beta_k^\th$ elements of $t$, then $s=j(f)(\tau)$. Consequently, $s\in
j(A)$ if and only if $j(f)(\tau)\in j(A)$, which holds if and only if $\tau\in j(f^{-1}A)$. This last property holds if and only if $f^{-1}A\in\mu$,
which can be computed in $V$. Therefore, I have shown that every $\delta$ approximation to $E$ over $V$ is in $V$, and so by the $\delta$
approximation property, $E\intersect V$ is in $V$.\QED

I stress again that counterexamples show that the corresponding ground model extender embedding $j_{E\intersect V}$ is not necessarily the same as
$j\restrict V$.

For the remainder of this section, I will show that forcing extensions obtained by forcing with a closure point at $\delta$ exhibit the
$\delta^\plus$ approximation and $\delta^\plus$ cover properties. Such closure point forcing extensions, therefore, fall under the scope of the Main
Theorem, and the results of this article consequently generalize \cite{Hamkins2001:GapForcing}. An abundance of reverse Easton iterations in the
literature, such as the Laver preparation or the Silver iteration to add Cohen subsets to regular cardinals, admit numerous closure points, and so
the Main Theorem is applicable. Recall that a poset $\Q$ is {\df $\leqdelta$ strategically closed} if there is a strategy for the second player in
the game of length $\delta+1$ allowing her to continue play, where the players alternate to build a descending sequence in $\Q$, with the second
player playing at limit stages. By {\df nontrivial} forcing, I mean one that necessarily adds a new set.

\Definition. A forcing notion has a {\df closure point} at $\delta$ when it factors as $\P*\Qdot$, where $\P$ is nontrivial, $\card{\P}\leq\delta$
and $\forces\Qdot$ is $\leqdelta$ strategically closed.

\Lemma. Forcing with a closure point at $\delta$ satisfies the $\delta^\plus$ approximation and $\delta^\plus$ cover
properties.\label{ApproximationLemma}

\Proof: Suppose that $V[g][H]$ has a closure point at $\delta$, so that $g*H\subseteq\P*\Qdot$ is $V$-generic, $\card{\P}\leq\delta$ and
$\Pforces``\Qdot$ is $\leqdelta$ strategically closed''. The $\delta^\plus$ cover property is easy to verify for $V\of V[g][H]$, because it holds
separately for each step of the forcing. For the $\delta^\plus$ approximation property, we reduce to the case of sets of ordinals, or binary ordinal
sequences, simply by enumerating sets in $V$ and considering approximations on the indices. So, suppose a sequence $s\in 2^\theta$ is in $V[g][H]$
and $s\restrict\sigma\in V$ whenever $\sigma$ has size at most $\delta$ in $V$. We want to show that $s$ itself is in $V$. By induction, we may
assume that all proper initial segments of $s$ are in $V$.

The easy case occurs when $\cof(\theta)\leq\delta$. It follows that $s\in V[g]$, and so there is a $\P$-name $\sdot$ in $V$ with $s=\sdot_g$. In $V$,
let $T$ be the tree of all possible initial segments of $s$, that is, $T=\set{t\in 2^\lttheta\st \boolval{\check t\of\sdot}^\P\neq 0}$. The sequence
$s$ is a branch through $T$ in $V[g]$. Since incomparable elements of this tree give rise to incompatible elements of $\P$, it is easy to see that
there are at most $\delta$ many {\df branch points} in $T$, elements $t\in T$ such that $t\concat0$ and $t\concat 1$ are both in $T$. Thus, the set
$\sigma$, consisting of the lengths of any such branch point, has size at most $\delta$ in $V$, and so $s\restrict\sigma$ is in $V$. But
$s\restrict\sigma$ gives exactly the information one needs to know, specifying which way to turn at any branch point, in order to follow the branch
$s$ through $T$. So $s\in V$, as desired.

For the remaining case, assume $\cof(\theta)>\delta$. Settling this case is exactly \cite[Key Lemma]{Hamkins2001:GapForcing}, but we include the
proof here for convenience. The idea is simply that if $s\notin V$, then any new small set $h$ added by $\P$ is forced by the closure of $\Qdot$ to
be embedded into $s$, giving an approximation not in the ground model. Let $\sdot$ be a $\P*\Qdot$-name for $s$, and suppose $\<p_0,\qdot_0>\in g*H$
forces that $\sdot$ is not in $\check V$, but all proper initial segments of $\sdot$ are in $\check V$. For each $\lambda<\theta$, choose in
$V[g][H]$ a condition $\<p_\lambda,\qdot_\lambda>\in g*H$ deciding $\sdot\restrict\lambda$ in $V$. Since $|\P|\leq\delta$, there must be a single
condition repeated for unboundedly many $p_\lambda$, and so we may in fact assume that $p_\lambda=p_0$ for all $\lambda$. By strengthening further if
necessary, we may assume that $\<p_0,\qdot_0>$ forces that $p_0$ has this property. Thus, for any $\lambda<\theta$ and any condition of the form
$\<p_0,\qdot>\leq\<p_0,\qdot_0>$, there is a stronger condition $\<p_0,\rdot>\leq\<p_0,\qdot>$ deciding $\sdot\restrict\check\lambda$. Since $\P$ is
nontrivial, there is some $h\in (2^\beta)^{V[g]}$ for some $\beta\leq\delta$ with $h\notin V$, but all initial segments of $h$ are in $V$. Let
$\dot\sigma$ be the name of a strategy witnessing that $\Qdot$ is $\leqdelta$ strategically closed. We construct in $V$ a tree of names $\qdot_t$ for
$t\in 2^\ltbeta$ for possible moves for the first player in the $\Qdot$ game, with the second player obeying the strategy $\dot\sigma$. Player one
begins with $\qdot_\emptyset=\qdot_0$. If $\qdot_t$ is defined, let $\rdot_t$ name the result of applying the strategy $\dot\sigma$ to the already
constructed play $\<\qdot_u\st u\of t>$ for player one, and let $b_t\in 2^\lttheta$ be the longest binary sequence such that
$\<p_0,\rdot_t>\forces\check b_t\of\sdot$. By our assumption on $\<p_0,\qdot_0>$, there are conditions $\qdot_{t\concat 0}$ and $\qdot_{t\concat 1}$
such that $\<p_0,\qdot_{t\concat i}>\leq\<p_0,\qdot_t>$ and $\<p_0,\qdot_{t\concat i}>\forces \check b_t\concat\check i\of\sdot$. If $t$ has limit
length and $\qdot_u$ is defined for all $u\ofnoteq t$, then  because these name conditions corresponding to a play according to $\dot\sigma$, there
is a name $\rdot_t$ for the result of applying $\dot\sigma$ to that play, and we let $\qdot_t$ name any stronger condition. In $V[g]$, after
interpreting the names, the sequence $\<q_t\st t\ofnoteq h>$ gives by construction the moves of player one in a play of length $\beta$ in $\Q$
according to $\sigma$, and so there is a condition $q\leq q_t$ for all $t\of h$. Thus, $q$ forces that $b=\union_{t\of h} b_t$ is a proper initial
segment of $s$, and so $b\in V$. But from $b$ we can reconstruct $h$ in $V$ by observing that $t$ deviates from $h$ exactly when $b_t$ deviates from
$b$. This contradicts our assumption that $h\notin V$.\QED

Mitchell \cite{Mitchell2003:ANoteOnHamkinsApproximationLemma} has provided a proof of (a generalization of) Lemma \ref{ApproximationLemma} that
avoids the tree construction by using master conditions.

\Corollary. The conclusions of the Main Theorem and its consequences hold for embeddings in any closure point forcing
extension.\label{ClosurePointForcing}

\Section Consequences of the Main Theorem\label{ConsequencesSection}

I will now apply the Main Theorem to a variety of large cardinal notions in order to show that if an extension satisfies the approximation and cover
properties, then it contains no new large cardinals. The case of the smaller large cardinals makes use of the following two lemmas.

\Lemma. Suppose that $V\of\Vbar$ satisfies the $\delta$ approximation and cover properties. If $\Xbar^\ltdelta\of\Xbar$ in $\Vbar$ and $\Xbar\elesub
\Vbar_\theta$ in the language with a predicate for $V$, so that $\<\Xbar,X,{\in}>\elesub\<\Vbar_\theta,V_\theta,{\in}>$, where $X=\Xbar\intersect V$,
then $X\in V$. Further, if $\Mbar$ is the Mostowski collapse of $\Xbar$, then the Mostowski collapse of $X$ is the same as $\Mbar\intersect
V$.\label{CondensationLemma}

\Proof: First, I will show $X\in V$. Suppose that $a\in V$ has size less than $\delta$ in $V$. Since $X\intersect a$ is a subset of $X$ of size less
than $\delta$ in $\Vbar$, it is in $\Xbar$. And since it is an element of $\Vbar_\theta$ of size less than $\delta$, it is covered by an element
$b\in V_\theta$ of size less than $\delta$ in $V_\theta$. By elementarity there is such a $b$ in $X$. Since $b$ has size less than $\delta$ and
$\delta\of X$, it follows that $b\of X$. In summary, we have $X\intersect a\of b\of X$, which implies $X\intersect a=b\intersect a$, and so
$X\intersect a$ is in $V$. Since all the $\delta$ approximations to $X$ over $V$ are in $V$, it follows by the $\delta$ approximation property that
$X\in V$.

Now consider $\<\Mbar,M,{\in}>$, the Mostowski collapse of $\<\Xbar,X,{\in}>$. Since $\Vbar_\theta$ knows that $V_\theta$ is transitive, it follows
that every element of $\Xbar$ that is an element of an element of $X$ is itself in $X$, and so the Mostowski collapse of $X$ is the same as the image
of $X$ under the Mostowski collapse of $\Xbar$; that is, $M$ is the Mostowski collapse of $X$. In particular, $M\in V$. It follows that $M\of
\Mbar\intersect V$. For the converse inclusion, let $\pi:\Xbar\iso\Mbar$ be the Mostowski collapse of $\Xbar$ and suppose that
$\pi(A)\in\Mbar\intersect V$, where $A\in\Xbar$. I may assume inductively that every element of $\pi(A)$ is in $M$. Thus, $A\intersect \Xbar\of X$.
It follows by elementarity that $A\of V$. Suppose that $a\in X$ has size less than $\delta$ in $X$. It follows that $A\intersect a$ is an element of
$\Xbar$, of size less than $\delta$ there, and a subset of $X$. Consequently, by the cover property, $A\intersect a\of b$ for some $b\in X$ of size
less than $\delta$ in $X$. Enumerate $b=\set{b_\alpha\st\alpha<\beta}$ in $V$, where $\beta<\delta$, and let $A_0=\set{\alpha\st b_\alpha\in
A\intersect a}$. Since $\pi$ fixes all ordinals below $\delta$ and all subsets of $\delta$, we see that $\alpha\in \pi(A_0)=A_0$ if and only if
$\pi(b_\alpha)\in \pi(A)\intersect\pi(a)$. Since these latter sets are all in $V$, including the sequence $\<\pi(b_\alpha)\st\alpha<\beta>$, it
follows that $A_0$ is in $V$, and consequently $A\intersect a\in V$. Thus also $A\intersect a\in X$, and so all $\delta$ approximations to $A$ using
$a\in X$ are in $X$. By elementarity, it follows that all $\delta$ approximations to $A$ over $V$ are in $V$, and so by the approximation property,
we conclude $A\in V$. This implies $\pi(A)\in M$, as desired.\QED

If the hypotheses concerning the approximation and cover properties are omitted from Lemma \ref{CondensationLemma}, then the conclusion can fail. For
example, if one adds a Prikry sequence $s$ to a measurable cardinal $\kappa>\delta$, then for any $\theta\geq\kappa$ there are elementary
substructures $\Xbar\elesub V_\theta[s]$ of size $\delta^\ltdelta$ in $V[s]$ such that $\Xbar^\ltdelta\of\Xbar$ in $V[s]$ and $s\in\Xbar$. In this
case, $\Xbar\intersect V$ is not in $V$, as it has size at most $\delta^\ltdelta$ but is unbounded in $\kappa$, violating the regularity of $\kappa$
in $V$.

\Lemma. Suppose that $V\of\Vbar$ satisfies the $\delta$ approximation and cover properties and $\kappa\geq\delta$ is an inaccessible cardinal. If
$A\of\kappa$ is any set in $\Vbar$, then there is a transitive model of set theory $\Mbar$ of size $\kappa$ in $\Vbar$ such that $A\in\Mbar$,
$\Mbar^\ltkappa\of\Mbar$ and $M=\Mbar\intersect V\in V$ is a model of set theory.\label{Mexists}

\Proof: Suppose $\ZFC^*$ is the fixed finite fragment of \ZFC\ used to define the models of set theory. The proof of the well-known L\'evy reflection
theorem establishes that there is an ordinal $\theta$ above $\kappa$ such that every formula appearing in $\ZFC^*$ reflects from the structure
$\<\Vbar,V,{\in}>$ to $\<\Vbar_\theta,V_\theta,{\in}>$. In particular, both $\Vbar_\theta$ and $V_\theta$ are models of set theory. In $\Vbar$ let
$\Xbar\elesub\Vbar_\theta$ be an elementary substructure of size $\kappa$ in the language with a predicate for $V$, so that
$\<\Xbar,X,{\in}>\elesub\<\Vbar_\theta,V_\theta,{\in}>$, where $X=\Xbar\intersect V$, such that $\Xbar^\ltkappa\of\Xbar$ and $A\in\Xbar$. By Lemma
\ref{CondensationLemma} the collapse $\Mbar$ of $\Xbar$ has the property that $M=\Mbar\intersect V$ is in $V$. And since $M$ is the collapse of $X$,
it is a model of set theory, as desired.\QED

\Corollary. Suppose $V\of\Vbar$ satisfies the $\delta$ approximation and cover properties. Then every weakly compact cardinal above $\delta$ in
$\Vbar$ is weakly compact in $V$.\label{CorollaryWC}

\Proof: Suppose $\kappa$ is weakly compact in $\Vbar$. For any subset $A\of\kappa$ in $V$ there is by Lemma \ref{Mexists} a model of set theory
$\Mbar$ in $\Vbar$ such that $A\in\Mbar$, $\Mbar^\ltkappa\of\Mbar$ and $M=\Mbar\intersect V$ is a model of set theory in $V$. Since $\kappa$ is
weakly compact in $\Vbar$, there is an embedding $j:\Mbar\to\Nbar$ in $\Vbar$ with critical point $\kappa$, and by using the induced normal
$\Mbar$-measure, we may assume $\Nbar^\ltkappa\of\Nbar$ in $\Vbar$. Since this embedding satisfies the hypotheses of the Main Theorem, it follows
that $j\restrict M:M\to N$ is an embedding in $V$. Since this restricted embedding still has critical point $\kappa$ and $A\in M$, it follows that
$\kappa$ is weakly compact in $V$.\QED

While the proof of the next theorem does not rely on the Main Theorem and the hypotheses are weaker, the result fits into the sequence of this
section. Recall that a cardinal $\kappa$ is {\df ineffable} if for any sequence $\<A_\alpha\st\alpha<\kappa>$ with $A_\alpha\of\alpha$ there is a set
$A\of\kappa$ such that $\set{\alpha<\kappa\st A_\alpha=A\intersect\alpha}$ is stationary.

\Theorem. Suppose $V\of\Vbar$ satisfies the $\delta$ approximation property. Then every ineffable cardinal $\kappa\geq\delta$ in $\Vbar$ is ineffable
in $V$.

\Proof: Suppose $\kappa\geq\delta$ is ineffable in $\Vbar$ and $\<A_\alpha\st\alpha<\kappa>$ is a sequence in $V$ with $A_\alpha\of\alpha$. In
$\Vbar$, there is a coherence set $A\of\kappa$ such that $B=\set{\alpha<\kappa\st A_\alpha=A\intersect\alpha}$ is stationary. In particular, all the
initial segments of $A$ are in $V$, and so by the $\delta$ approximation property, the set $A$ itself is in $V$. It follows that $B\in V$ also, and
there can be no club in $V$ avoiding $B$, as there is no such club in $\Vbar$. So $\kappa$ is ineffable in $V$.\QED

Kai Hauser \cite{Hauser1991:IndescribableCardinals} provided a useful embedding characterization of indescribable cardinals, to which the Main
Theorem applies, by showing for natural numbers $m,n\geq 1$ that a cardinal $\kappa$ is {\df $\Pi^m_n$ indescribable} if for any transitive model of
set theory $M$ of size $\kappa$ with $M^\ltkappa\of M$ and $\kappa\in M$, there is a transitive set $N$ and an elementary embedding $j:M\to N$ with
critical point $\kappa$ such that $N$ is $\Sigma^m_{n-1}$ correct, that is, such that $(V_{\kappa+m})^N\elesub_{n-1} V_{\kappa+m}$ and
$N^\card{V_{\kappa+m-2}}\of N$ (meaning $N^\ltkappa\of N$ when $m=1$). Since any first order statement about $V_{\kappa+m}$ is $\Delta_0$ in
$V_{\kappa+m+1}$, using $V_{\kappa+m}$ as a parameter, it follows that $\Pi^{m+1}_1$ indescribability implies $\Pi^m_n$ indescribability for any $n$.
A cardinal $\kappa$ is {\df totally indescribable} when it is $\Pi^m_n$ indescribable for any $m,n\in\omega$, or equivalently, if it is $\Pi^m_1$
indescribable for every $m$. Since $\Pi^1_1$ indescribability is simply weak compactness, the next corollary generalizes Corollary \ref{CorollaryWC}.

\Corollary. Suppose $V\of\Vbar$ satisfies the $\delta$ approximation and cover properties. Then every totally indescribable cardinal above $\delta$
in $\Vbar$ is totally indescribable in $V$. Indeed, for $m\geq 1$ every $\Pi^m_1$ indescribable cardinal above $\delta$ in $\Vbar$ is $\Pi^m_1$
indescribable in $V$.

\Proof: Suppose that $\kappa$ is $\Pi^m_1$ indescribable in $\Vbar$, and consider any transitive model of set theory $M_0$ in $V$ with
$M_0^\ltkappa\of M_0$ and $\kappa\in M_0$. By Lemma \ref{Mexists}, there is a transitive model of set theory $\Mbar$ in $\Vbar$ with
$\Mbar^\ltkappa\of\Mbar$ in $\Vbar$ and $M_0\in\Mbar$ such that $M=\Mbar\intersect V$ is also a model of set theory. Since $\kappa$ is $\Pi^m_1$
indescribable in $\Vbar$, there is an embedding $j:\Mbar\to\Nbar$ such that $\Nbar$ is $\Sigma^m_0$ correct in $\Vbar$. By the Main Theorem, the
restricted embedding $j\restrict M:M\to N$ lies in $V$. By restricting the embedding further, down to $M_0$, I obtain the embedding $j_0=j\restrict
M_0:M_0\to N_0$, where $N_0=j(M_0)$.

I claim that $N_0$ is $\Sigma^m_0$ correct. To see this, observe first by Corollary \ref{Extra} that $N^{|V_{\kappa+m-2}|}\of N$ in $V$, since
$V_{\kappa+m-2}\of\Vbar_{\kappa+m-2}$ and $\Nbar^{|\Vbar_{\kappa+m-2}|}\of\Nbar$ in $\Vbar$, and consequently
$|V_{\kappa+m-2}|^V<(|\Vbar_{\kappa+m-2}|^\plus)^\Vbar$. Since $M$ knows that $M_0^\ltkappa\of M_0$, it follows that $N$ knows that
$N_0^{<j(\kappa)}\of N_0$. Because $N$ has all the sequences over $N_0$ of length up to $|V_{\kappa+m-2}|$, which is less than $j(\kappa)$, it
follows that $N_0^{|V_{\kappa+m-2}|}\of N_0$ in $V$, as required. Second, because $(N_0)_{\kappa+m}$ is a transitive subset of $V_{\kappa+m}$, it
follows that $\Sigma_0$ truth is preserved. So the embedding $j_0:M_0\to N_0$ is $\Sigma^m_0$ correct, and the proof is complete.\QED

Recall from \cite{Villaveces98} that a cardinal $\kappa$ is {\df unfoldable} if it is {\df $\theta$ unfoldable} for every ordinal $\theta$, meaning
that for any transitive model of set theory $M$ of size $\kappa$ there is a transitive set $N$ and an embedding $j:M\to N$ with critical point
$\kappa$ such that $j(\kappa)>\theta$. It suffices if such embeddings $j$ exist for arbitrarily large sets $M$, that is, if every $A\of\kappa$ can be
placed into such an $M$ (proof: given any $M'$, place it into an $M$, get the embedding and restrict it to $M'$). The cardinal $\kappa$ is {\df
strongly unfoldable} if it is $\theta$ strongly unfoldable for every ordinal $\theta$, meaning that for every transitive model of set theory $M$ of
size $\kappa$ with $M^\ltkappa\of M$ there is an embedding $j:M\to N$ into a transitive set $N$ with critical point $\kappa$ such that
$j(\kappa)>\theta$ and $V_\theta\of N$. If $\theta$ is a successor ordinal or has cofinality above $\kappa$, such an $N$ can be found for which
$N^\kappa\of N$ (see \cite{Hamkins:LaverDiamond}).

\Corollary. Suppose $V\of\Vbar$ satisfies the $\delta$ approximation and cover properties. Then every strongly unfoldable cardinal above $\delta$ in
$\Vbar$ is strongly unfoldable in $V$.

\Proof: Fix any successor ordinal $\theta$ and any set $A\of\kappa$ in $V$. By Lemma \ref{Mexists} there is a transitive model $\Mbar$ of size
$\kappa$ in $\Vbar$ such that $\Mbar^\ltkappa\of\Mbar$ in $\Vbar$, $A\in\Mbar$ and $M=\Mbar\intersect V$ is a model of set theory in $V$. Since
$\kappa$ is $\theta$ strongly unfoldable in $\Vbar$, there is an embedding $j:\Mbar\to\Nbar$ with $\Vbar_\theta\of\Nbar$ and $\Nbar^\kappa\of\Nbar$
in $\Vbar$. Thus, the Main Theorem applies, and so the restricted embedding $j\restrict M:M\to N$ exists in $V$. By Corollary \ref{Extra} it follows
that $V_\theta\of N$, and we know $A\in M$, so this restricted embedding serves to witness the $\theta$ strong unfoldability (for $A$) in $V$.\QED

We don't actually need $\Nbar^\kappa\of\Nbar$ in the previous argument, but rather only $\Nbar^\delta\of\Nbar$. And since such an embedding can be
found when $\theta$ is either a successor ordinal or has cofinality above $\delta$, we conclude the following.

\Corollary. Suppose $V\of\Vbar$ satisfies the $\delta$ approximation and cover properties and $\theta$ is a successor ordinal or has cofinality above
$\delta$. Then every  $\theta$ strongly unfoldable cardinal above $\delta$ in $\Vbar$ is $\theta$ strongly unfoldable in $V$.

\Corollary. Suppose $V\of\Vbar$ satisfies the $\delta$ approximation and cover properties. Then every measurable cardinal above $\delta$ in $\Vbar$
is measurable in $V$.

\Proof: Since $\kappa$ is measurable in $\Vbar$, there is a normal ultrapower embedding $j:\Vbar\to\Nbar$. Since $\Nbar^\kappa\of\Nbar$, the Main
Theorem implies that the restricted embedding $j\restrict V:V\to N$ is amenable to $V$. In $V$ one may construct a normal measure $\mu$ on $\kappa$
from $j\restrict P(\kappa)^V$ by defining $X\in\mu\iff\kappa\in j(X)$.\QED

A cardinal $\kappa$ is {\df tall} if it is $\theta$ tall for every $\theta$, meaning that there is an embedding $j:V\to M$ with critical point
$\kappa$ such that $j(\kappa)>\theta$ and $M^\kappa\of M$.

\Corollary. Suppose $V\of\Vbar$ satisfies the $\delta$ approximation and cover properties. Then every tall cardinal above $\delta$ in $\Vbar$ is tall
in $V$. Indeed, for any $\theta$, every $\theta$ tall cardinal above $\delta$ in $\Vbar$ is $\theta$ tall in $V$.

\Proof: This is immediate when we have the $\delta$ approximation property for classes, because if a class $j:\Vbar\to\Nbar$ witnesses that $\kappa$
is $\theta$ tall in $\Vbar$, then the restriction $j\restrict V:V\to N$ witnesses that $\kappa$ is $\theta$ tall in $V$. But in general we only have
amenability, so we work with the induced extenders. Suppose that $j:\Vbar\to\Nbar$ witnesses that $\kappa$ is $\theta$ tall in $\Vbar$. By the Main
Theorem, we know that $j\restrict V:V\to N$ is amenable to $V$, where $N=\Nbar\intersect V$. By Corollary \ref{Extra}, we know $N^\kappa\intersect
V\of N$. Let $E=j\restrict P(\kappa)^V$, which is in $V$ by amenability, and let $j_E:V\to N_E$ be the corresponding extender embedding. Thus,
$j_E\restrict P(\kappa)^V=j\restrict P(\kappa)^V$ and every element of $N_E$ has the form $j_E(f)(\beta)$ for some $f\in V^\kappa\intersect V$ and
$\beta<j(\kappa)$. Now suppose that $\<j_E(f_\alpha)(\beta_\alpha)\st\alpha<\kappa>$ is a $\kappa$ sequence from $N_E$ in $V$. We may assume
$\<\beta_\alpha\st\alpha<\kappa>$ is also in $V$, and so by the $\kappa$ closure of $\Nbar$ it is in $\Nbar\intersect V=N$. Because $N_E$ and $N$
agree up to rank $j(\kappa)$, this means that $\<\beta_\alpha\st\alpha<\kappa>\in N_E$. Since
$\<j_E(f_\alpha)\st\alpha<\kappa>=j(\<f_\alpha\st\alpha<\kappa>)\restrict\kappa$ is in $N_E$ as well, we see that
$\<j(f_\alpha)(\beta_\alpha)\st\alpha<\kappa>$ is in $N_E$, as desired.\QED

\Corollary. Suppose $V\of\Vbar$ satisfies the $\delta$ approximation and cover properties. Then every strong cardinal above $\delta$ in $\Vbar$ is
strong in $V$.

\Proof: Suppose that $\kappa$ is $\theta$ strong in $\Vbar$ and $\theta$ is either a successor ordinal or has cofinality above $\delta$. In $\Vbar$
there is a $\theta$ strongness extender embedding $j:\Vbar\to\Nbar$ with $\cp(j)=\kappa$, $\Vbar_\theta\of\Nbar$ and $\Nbar^\delta\of\Nbar$. By the
Main Theorem, the restricted embedding $j\restrict V:V\to N$ is amenable to $V$, and by Corollary \ref{Extra}, we know $V_\theta\of N$. Let
$E=j\restrict P(\kappa)^V$ in $V$ and observe that the corresponding extender embedding $j_E:V\to M_E$ has $j_E\restrict P(\kappa)^V=j\restrict
P(\kappa)^V$, and consequently $V_\theta\of M_E$. So $\kappa$ is $\theta$ strong in $V$.\QED

\Corollary. Suppose $V\of\Vbar$ satisfies the $\delta$ approximation and cover properties. Then every Woodin cardinal above $\delta$ in $\Vbar$ is
Woodin in $V$.

\Proof: If $\kappa$ is Woodin in $\Vbar$, then for every $A\of\kappa$ in $V$ there is $\gamma\in(\delta,\kappa)$ such that for arbitrarily large
$\lambda<\kappa$ there is an extender embedding $j:\Vbar\to\Nbar$ such that $\cp(j)=\gamma$ and $j(A)\intersect\lambda=A\intersect\lambda$. We may
also assume $\Nbar^\gamma\of\Nbar$. It follows from the Main Theorem that the restriction $j\restrict V:V\to N$ is amenable to $V$. And of course it
still satisfies $j(A)\intersect\lambda=A\intersect\lambda$. Since $j\restrict P(\kappa)^V\in V$ by amenability, the induced extender embeddings
therefore witness that $\kappa$ is a Woodin cardinal in $V$.\QED

The case of strongly compact cardinals presents peculiar difficulties, and it will be treated separately in Section \ref{StrCSection}. So I move now
to the case of supercompact cardinals.

\Corollary. Suppose $V\of\Vbar$ satisfies the $\delta$ approximation and cover properties. Then every supercompact cardinal above $\delta$ in $\Vbar$
is supercompact in $V$. Indeed, for any $\theta$, every $\theta$ supercompact cardinal above $\delta$ in $\Vbar$ is $\theta$ supercompact in
$V$.\label{SupercompactCorollary}

\Proof: If $j:\Vbar\to\Nbar$ is an embedding in $\Vbar$ witnessing that $\kappa$ is $\theta$ supercompact in $\Vbar$, then by the Main Theorem the
restriction $j\restrict V:V\to N$ is amenable to $V$ and $N=\Nbar\intersect V$. By Corollary \ref{Extra} we know $N^\theta\of N$ in $V$, and so
$j\image\theta\in N$. Thus, from $j\restrict P(P_\kappa\theta)^V$ we may in $V$ construct the induced normal fine measure $\mu$ on $P_\kappa\theta$
by defining $X\in \mu\iff j\image\theta\in j(X)$. So $\kappa$ is $\theta$ supercompact in $V$.\QED

\Corollary. Suppose $V\of\Vbar$ satisfies the $\delta$ approximation and cover properties. Then every almost huge, huge or superhuge cardinal above
$\delta$ in $\Vbar$ exhibits the same large cardinal property in $V$.

\Proof: Once again, suitable restrictions of these embeddings witness the large cardinal property in $V$.\QED

Let me close this section with some results on the question of making a weakly compact or measurable cardinal $\kappa$ indestructible by
$\ltkappa$-directed closed forcing. The only method currently known for doing this is to begin with a supercompact cardinal $\kappa$ and perform the
Laver preparation. The following theorem shows that if one produces an indestructible weakly compact cardinal in any extension resembling the Laver
preparation, that is, one exhibiting the approximation and cover properties, then one must have begun with a supercompact cardinal. This theorem
generalizes a result in \cite{ApterHamkins2001:IndestructibleWC}.

\Theorem. Suppose $V\of\Vbar$ satisfies the $\delta$ approximation and cover properties for some $\delta<\kappa$. If $\kappa$ is weakly compact in
$\Vbar$ and $(2^{\theta^\ltkappa})^V$ is collapsed to $\kappa$ in $\Vbar$, then $\kappa$ was $\theta$ supercompact in
$V$.\label{IndestructibleWCtheta}

\Proof: By Lemma \ref{Mexists} there is a transitive model of set theory $\Mbar$ in $\Vbar$ of size $\kappa$ such that $P(P_\kappa\theta)^V\in\Mbar$,
$\Mbar^\ltkappa\of\Mbar$, $M=\Mbar\intersect V$ is a model of set theory and $\Mbar$ knows that $|\theta|=\kappa$. Since $\kappa$ is weakly compact
in $\Vbar$ there is an embedding $j:\Mbar\to \Nbar$ with critical point $\kappa$ and $\Nbar^\ltkappa\of\Nbar$. Since $\theta$ has size $\kappa$ in
$\Mbar$ there is a relation $\trianglelt$ on $\kappa$ with order type $\theta$. Notice that if $\beta<\kappa$ has order type $\alpha$ with respect to
$\trianglelt$, then $j(\beta)=\beta$ has order type $j(\alpha)$ with respect to $j(\trianglelt)$. Therefore, $j\image\theta$ is constructible in
$\Nbar$ from $\trianglelt$ and $j(\trianglelt)$, and so $j\image\theta\in \Nbar$. By the Main Theorem, $j\restrict M:M\to N$, where
$N=\Nbar\intersect V$, is an embedding in $V$. In particular, $j\image\theta$ is in $V$, and hence in $N$. In $V$, the set $\mu$ of all $X\of
P_\kappa\theta$ such that $j\image\theta\in j(X)$ is a normal fine measure on $P_\kappa\theta$, and so $\kappa$ is $\theta$ supercompact there.\QED

Using the results on closure point forcing, we obtain the following corollary, one of the main theorems of \cite{ApterHamkins2001:IndestructibleWC}.

\Corollary.(\cite{ApterHamkins2001:IndestructibleWC}) If $V\of V[G]$ has a closure point below $\kappa$ and the weak compactness of $\kappa$ is
indestructible over $V[G]$ by the forcing to collapse cardinals to $\kappa$, then $\kappa$ was supercompact in $V$.

\Proof: Such extensions, when followed by further $\ltkappa$ directed closed forcing, still have the same closure point at some $\delta<\kappa$, and
consequently by Lemma \ref{ApproximationLemma} exhibit the $\delta^\plus$ approximation and cover properties. So the corollary follows from the
previous theorem.\QED

The same idea applies to indestructible measurable cardinals, but here one only needs to know that $\theta$ is collapsed to $\kappa$, rather than
$2^{\theta^\ltkappa}$ as above.

\Theorem. Suppose $V\of\Vbar$ satisfies the $\delta$ approximation and cover properties. If $\kappa>\delta$ is measurable in $\Vbar$ and $\theta$ has
cardinality $\kappa$ in $\Vbar$, then $\kappa$ was $\theta$ supercompact in $V$.\label{ThetaSCinV}

\Proof: Let $j:\Vbar\to\Nbar$ be the ultrapower embedding by a normal measure on $\kappa$ in $\Vbar$. It follows by the Main Theorem that the
restriction $j\restrict V:V\to N$ is amenable to $V$, and so $j\image\theta\in V$. Furthermore, since $|\theta|=\kappa$ in $\Vbar$ and
$\Nbar^\kappa\of\Nbar$ in $\Vbar$, it follows that $j\image\theta\in\Nbar$, and so $j\image\theta\in \Nbar\intersect V=N$. From $j\restrict
P(P_\kappa\theta)^V$ in $V$ one can therefore construct a normal fine measure $\mu$ by defining $X\in\mu\iff j\image\theta\in j(X)$. So $\kappa$ is
$\theta$ supercompact in $V$.\QED

\Section The Case of Strongly Compact Cardinals\label{StrCSection}

The case of strongly compact cardinals presents special problems for the arguments of Section \ref{ConsequencesSection}, the main obstacle being that
the restriction $j\restrict V$ of a strong compactness embedding $j:\Vbar\to\Nbar$ does not seem immediately to reveal the full strength of the
original embedding, as it did so easily in the case of measurability, supercompactness and so on. Here, in order to carry out the analysis for
strongly compact cardinals, I will make some additional assumptions about the nature of the extension $V\of\Vbar$.

\Theorem. Suppose $\delta<\kappa$ and $V\of\Vbar$ satisfies the $\delta$ approximation and cover properties, as well as the $\kappa$ cover property.
For any $\theta$, if $\kappa$ is $\theta$ strongly compact in $\Vbar$, then it was $\theta$ strongly compact in $V$.\label{KappaCovering}

\Proof: Suppose that $j:\Vbar\to\Nbar$ is a $\theta$ strong compactness embedding in $\Vbar$, the ultrapower by a fine measure $\mu$ on
$P_\kappa\theta$. Let $s=[\id]_\mu$, so that $j\image\theta\of s\of j(\theta)$ and $|s|<j(\kappa)$ in $\Nbar$. Since $\Nbar^\kappa\of\Nbar$ in
$\Vbar$, the Main Theorem establishes that $j\restrict V:V\to N$, where $N=\Nbar\intersect V$, is amenable to $V$. By $j$ of the $\kappa$ cover
property, it follows that $s\of t$ for some $t\in N$ of size less than $j(\kappa)$ in $N$. Without loss of generality, $t\of j(\theta)$. Since also
$j\image\theta\of t$, it follows that $t$ generates a fine measure $\mu$ on $P_\kappa\theta$ in $V$, defined by $X\in\mu\iff t\in j(X)$. So $\kappa$
is $\theta$ strongly compact in $V$.\QED

The $\kappa$ cover property of $V\of\Vbar$ captures the operative power of {\df mildness} in \cite[Corollary 16]{Hamkins2001:GapForcing}, where a
poset $\P$ is mild relative to $\kappa$ if every set of ordinals of size less than $\kappa$ in $V^\P$ has a nice name of size less than $\kappa$.
(The definition in \cite{Hamkins2001:GapForcing} erroneously omitted the requirement that the name be nice, though this was used in the proofs.) I
conjecture that the added assumption of $\kappa$ covering in Theorem \ref{KappaCovering} is unnecessary.

\Conjecture. Suppose $\delta<\kappa$ and $V\of\Vbar$ satisfies the $\delta$ approximation and cover properties. For any $\theta$, if $\kappa$ is
$\theta$ strongly compact in $\Vbar$, then it is $\theta$ strongly compact in $V$.\label{ConjectureOnStrC}

The results of this section point towards a positive resolution of this conjecture. First, I can improve Theorem \ref{KappaCovering} by weakening the
assumption of $\kappa$ covering to the assumption only that no regular cardinal above $\kappa$ in $V$ has cofinality below $\kappa$ in $\Vbar$. This
argument will rely on an old characterization of strong compactness due to Ketonen. A filter on $\lambda$ is {\df uniform} if it contains the tail
segments $[\beta,\lambda)$ for every $\beta<\lambda$.

\Fact.({Ketonen \cite{Ketonen72:CardinalSins}}) A cardinal $\kappa$ is $\theta$ strongly compact if and only if for every regular cardinal
$\lambda\in[\kappa,\theta]$ there is a $\kappa$-complete uniform ultrafilter on $\lambda$.\label{KetonenFact}

Imagine, for example, that we have an embedding $j:V\to N$ with critical point $\kappa$ that is discontinuous at $\lambda$ in the sense that $\sup
j\image\lambda<j(\lambda)$. For any $\alpha\in[\sup j\image\lambda,j(\lambda))$ one may define a measure $\mu$ on $\lambda$ by $X\in\mu$ if and only
if $\alpha\in j(X)$, and it is easy to see that this will be a $\kappa$ complete uniform ultrafilter on $\lambda$. Conversely, the ultrapower by any
such measure $\mu$ will be discontinuous at $\lambda$, as $\sup j_\mu\image\lambda\leq[\id]_\mu<j_\mu(\lambda)$.

\Theorem. Suppose $\delta<\kappa\leq\theta$ and $V\of\Vbar$ satisfies the $\delta$ approximation and cover properties and every regular cardinal of
$V$ in the interval $(\kappa,\theta]$ has cofinality at least $\kappa$ in $\Vbar$. If $\kappa$ is $\theta$ strongly compact in $\Vbar$, then $\kappa$
was $\theta$ strongly compact in $V$.\label{StrC}

\Proof: Suppose that $\kappa$ is $\theta$ strongly compact in $\Vbar$. Fix a $\theta$ strong compactness ultrapower embedding $j:\Vbar\to\Nbar$ by a
fine measure $\mu$ on $P_\kappa\theta$. Let $s=[\id]_\mu$, so that $j\image\theta\of s\of j(\theta)$ and $|s|<j(\kappa)$ in $\Nbar$. Suppose
$\lambda$ is in the interval $[\kappa,\theta]$ and regular in $V$. By assumption, $\kappa\leq\cof(\lambda)$ in $\Vbar$. It follows that
$t=s\intersect j(\lambda)$, which has size less than $j(\kappa)$ in $\Nbar$, is bounded in $j(\lambda)$, and yet $j\image\lambda\of t$. Therefore
$\sup j\image\lambda<j(\lambda)$, and so $j$ is discontinuous at $\lambda$.

Since $j$ is the ultrapower by a measure on some set, it follows that $\Nbar^\kappa\of\Nbar$, and so the Main Theorem applies. Consequently, the
restricted embedding $j\restrict V:V\to N$, where $N=\Nbar\intersect V$, is amenable to $V$. Since the restricted embedding of course still satisfies
$\sup j\image\lambda<j(\lambda)$, one can use $j\restrict P(\lambda)^V$ as above to construct a $\kappa$-complete uniform ultrafilter on $\lambda$ in
$V$. By Ketonen's result, it follows that $\kappa$ is $\theta$ strongly compact in $V$.\QED

In particular, if $\Vbar$ preserves all cardinals and cofinalities over $V$, then the hypotheses of Theorem \ref{StrC} are satisfied, and so
Conjecture \ref{ConjectureOnStrC} holds for such extensions. Note that Theorem \ref{KappaCovering} is an immediate corollary of Theorem \ref{StrC},
because the $\kappa$ cover property implies that every cardinal with cofinality below $\kappa$ in $\Vbar$ has cofinality below $\kappa$ in $V$.

A perusal of Ketonen's argument \cite{Ketonen72:CardinalSins} will reveal that in order to establish Fact \ref{KetonenFact} one does not need that
{\it every} regular $\lambda$ has a $\kappa$ complete uniform ultrafilter on $\lambda$, but rather only that $\mu$-almost every $\lambda$ has that
property, where $\mu$ is any $\kappa$-complete uniform weakly normal ultrafilter on $\theta$, concentrating on cardinals of cofinality at least
$\kappa$. We may consequently also weaken the corresponding hypothesis of Theorem \ref{StrC}.

I would like next to observe that the critical exception making Theorem \ref{StrC} weaker than Conjecture \ref{ConjectureOnStrC}---the case of a
cardinal $\kappa$ that is $\theta$ strongly compact for a cardinal $\theta$ that was regular in $V$ but has cofinality less than $\kappa$ in
$\Vbar$---simply does not not occur with supercompactness. The situation here is rather similar to the fact that Prikry forcing above a strongly
compact cardinal destroys it. If one could extend Observation \ref{PrikryLikeObservation} to the case of strong compactness, this would prove that
Conjecture \ref{ConjectureOnStrC} is true.

\Observation. Suppose $V\of\Vbar$ satisfies the $\delta$ approximation and cover properties. If $\delta<\kappa\leq\theta$ and $\theta$ is a regular
cardinal of $V$ that has cofinality less than $\kappa$ in $\Vbar$, then $\kappa$ is not $\theta$ supercompact in
$\Vbar$.\label{PrikryLikeObservation}

\Proof: Suppose $\kappa$ is $\theta$ supercompact in $\Vbar$, so that there is a $\theta$ supercompactness embedding $j:\Vbar\to\Nbar$ in $\Vbar$. In
particular, $j\image\theta\in\Nbar$. Furthermore, since $\cof(\theta)<\kappa$, it follows that $\sup j\image\theta=j(\theta)$. By the Main Theorem,
however, the restricted embedding $j\restrict V:V\to N$ is amenable to $V$, and so $j\image\theta$ is in $V$. Consequently,
$j\image\theta\in\Nbar\intersect V=N$. Since $j\image\theta$ has size $\theta<j(\kappa)\leq j(\theta)$ and $j(\theta)$ is regular in $N$, it follows
that $\sup j\image\theta<j(\theta)$, a contradiction.\QED

Let me now prove the conjecture outright in the case of $\theta=\kappa^\plus$.

\Theorem. Suppose $\delta<\kappa$ and $V\of\Vbar$ satisfies the $\delta$ approximation and cover properties. If $\kappa$ is $\kappa^\plus$ strongly
compact in $\Vbar$, then it is $\kappa^\plus$ strongly compact in $V$ (interpreting $\kappa^\plus$ separately in $\Vbar$ and $V$,
respectively).\label{kappaplus}

\Proof: The essential idea here was employed in \cite[Lemma 2.3]{Apter2003:IndestructibilityStrongnessAndLevelByLevelEquivalence}. There are two
cases. If $\kappa^\plus$ is preserved from $V$ to $\Vbar$, this theorem is a special case of Theorem \ref{StrC}. Alternatively, if $\kappa^\plus$ is
collapsed from $V$ to $\Vbar$, then $\kappa$ is $\kappa^\plus$ supercompact in $V$ by Theorem \ref{ThetaSCinV}, and hence $\kappa^\plus$ strongly
compact there, as desired.\QED

Theorem \ref{kappaplus} does not seem to rule out the possibility that the degree of strong compactness of $\kappa$ increased from $V$ to $\Vbar$,
since it appears to be compatible with the conclusion of the theorem that $\kappa$ is $(\kappa^\plus)^V$ strongly compact in $V$ (but not more) and
$(\kappa^\plus)^\Vbar$ strongly compact in $\Vbar$, even when $(\kappa^\plus)^V<(\kappa^\plus)^\Vbar$. Such a phenomenon, however, is exactly what
Conjecture \ref{ConjectureOnStrC} rules out.

So altogether the evidence in favor of Conjecture \ref{ConjectureOnStrC} is first, that it follows the pattern of the results for all the other large
cardinals in Section \ref{ConsequencesSection}, such as Corollary \ref{SupercompactCorollary}; second, it holds when there is a small additional
degree of cofinality preservation by Theorems \ref{KappaCovering} and \ref{StrC}; third, failures of this additional preservation are incompatible
with $\kappa$ being $\theta$ supercompact in $\Vbar$ by Observation \ref{PrikryLikeObservation}, suggesting that they might also be incompatible with
$\kappa$ being $\theta$ strongly compact; and finally, fourth, it holds outright in the case of $\theta=\kappa^\plus$ by Theorem \ref{kappaplus}.

I close the article with an application of the Main Theorem by showing that it provides a new, easier proof of the second main Theorem of
\cite{HamkinsShelah98:Dual}, improving it to the case of strategically closed forcing.

\Theorem. After any nontrivial forcing of size less than $\kappa$, any further $\ltkappa$ strategically closed forcing that adds a new subset to any
$\lambda$ will destroy the $\lambda$ strong compactness of $\kappa$.

\Proof: Suppose that $g*G\of\P*\Q$ is $V$ generic for forcing with $|\P|<\kappa$ and $\forces_\P\Qdot$ is $\ltkappa$ strategically closed. Let
$A\of\lambda$ be a set that is in $V[g][G]$ but not in $V[g]$. Suppose that $\kappa$ is $\lambda$ strongly compact in $V[g][G]$, so that there is an
ultrapower embedding $j:V[g][G]\to N[g][j(G)]$ by a fine measure $\mu$ on $P_\kappa\lambda$. By Lemma \ref{ApproximationLemma}, this forcing has the
$\delta$ approximation and cover properties, where $\delta=|\P|^\plus$. Further, the model $N[g][j(G)]$ is closed under $\kappa$ sequences in
$V[g][G]$, since $j$ is the ultrapower by a $\kappa$-complete measure. Thus, by the Main Theorem, $N\of V$ and $j\restrict V:V\to N$ is amenable to
$V$.

Let $s=[\id]_\mu$, so that $j\image\lambda\of s\of j(\lambda)$ and $|s|^{N[g][j(G)]}<j(\kappa)$. Since $j(\Q)$ is ${<}j(\kappa)$ strategically
closed, it follows that $s\in N[g]$, a small forcing extension, and so $s\of t$ for some $t\in N$ with $t\of j(\lambda)$ and $|t|^N<j(\kappa)$. Using
that $j\image\lambda\of t$, it follows that $\alpha\in A\iff j(\alpha)\in j(A)\iff j(\alpha)\in j(A)\intersect t$. And since $j(A)\intersect t$ is a
set of ordinals in $N[g][j(G)]$ of size less than $j(\kappa)$, it follows by the strategic closure of $j(\Q)$ that it is in $N[g]$, which is a
subclass of $V[g]$. Therefore, we may construct $A$ in $V[g]$ by $\alpha\in A\iff j(\alpha)\in j(A)\intersect t$, using the fact that
$j\restrict\lambda\in V$. So $A$ is in $V[g]$ after all, a contradiction.\QED

It follows that small forcing always kills Laver indestructibility. The theorem can be improved with the observation that the proof used only the
fact that $\Q$ was ${\leq}|\P|$ strategically closed and didn't add any new sequences of ordinal of length less than $\kappa$. This establishes:

\Theorem. After nontrivial forcing $\P$ of size $\delta<\kappa$, any further forcing $\Q$ which is $\leqdelta$ strategically closed and $\ltkappa$
distributive which adds a subset to any $\lambda$ destroys the $\lambda$ strong compactness of $\kappa$.

For example, if one adds a Cohen subset to $\delta$ and then to $\lambda$, one destroys all strongly compact cardinals in the interval
$(\delta,\lambda]$.

\bigskip\noindent
{\sc
The Graduate Center of The City University of New York\\
Mathematics Program, 365 Fifth Avenue, New York, NY 10016}\\
jdh@hamkins.org, http://jdh.hamkins.org\\

\bibliographystyle{alpha}
\bibliography{MathBiblio,HamkinsBiblio}

\end{document}